\documentclass{article}

\def \N {{\mathbf {N}}}

\def \G {{\mathbf {G}}}
\def \P {{\mathbf {P}}}
\def \B {{\mathcal B}}

\def\uu{\sqcup}
\title{ Asymmetric Mixing Poisson Suspensions}
\author{ Valery V. Ryzhikov}
\date{}

\begin{document}
\Large

\maketitle

\begin{abstract} There are mixing Poisson suspensions that are not isomorphic to their inverses.

\vspace{2mm}

\end{abstract}

\section{Poisson Suspensions}
Gaussian automorphisms and Poisson suspensions are well known in ergodic theory \cite{KSF} and other areas of mathematics \cite{Ne}.
The former are obtained by embedding $ \G$ of the orthogonal group into the group $Aut$ of automorphisms of a probability space. The latter arises from the embedding $ \P$ of the group $Aut_\infty$ into $Aut$, where $Aut_\infty$ is the group of (infinite) automorphisms of $T$ of the standard Lebesgue space with $(X,\B,\mu)$ and sigma-finite measure $\mu$. The automorphism $T$ induces an orthogonal operator on $L_2(\mu)$, which corresponds to a Gaussian automorphism $\G(T)$, and the Poisson suspension $\P(T)$ acts on the Poisson probability space.

The spectral properties of the Gaussian automorphism $\G(T)$ and the suspension $\P(T)$ are known to be identical (see \cite{Ne}). However, their algebraic properties can be different. A Gaussian automorphism $\G(T)$ is always conjugate to its inverse $\G(T^{-1})$, but the Poisson suspensions $\P(T)$ and $\P(T^{-1})$ may inherit an asymmetry in the case of 
non-isomorphism of $T$ and $T^{-1}$. Several contrasts between Gaussian and Poisson suspensions  are shown in \cite{PR}. There are  prime  Poisson suspensions (possessing  only trivial invariant sigma-subalgebras) with  trivial centralizer. Gaussian automorphism  is in a continuum of flows and have a very reach  structure of invariant sigma-subalgebras. Note that in \cite{R25}, Poisson suspensions without roots but with large centralizer is   presented.

Now we  briefly recall the definition of a Poisson suspension.
Consider the space $(X,\mu)$ with sigma-finite Lebesgue measure, where $X$ is the union of disjoint intervals. The configuration space $X_\circ$ is, by definition, formed by all infinite countable sets $x_\circ$ under the condition that each interval of $X$ contains only finitely many elements of the set $x_\circ$. The space $X_\circ$ is equipped with the Poisson measure $\mu_\circ$ as follows. Subsets $A\subset X$ of finite measure are associated with cylindrical sets
$C(A,k)=\{x_\circ\in X_\circ \ : \ |x_\circ\cap A|=k\},$ $k=0,1,2,\dots$.
For measurable disjoint sets $A_1, A_2,\dots, A_N$, we have
$$\mu_\circ(\bigcap_{i=1}^N C(A_i,k_i))=\prod_{i=1}^N \frac {\mu(A_i)^{k_i}}{k_i!} e ^{-\mu(A_i)}.$$ Continuing the measure $\mu_\circ$ from the semiring of cylinders leads to the Poisson space $(X_\circ,\mu_\circ)$.

An automorphism $T$ of the space $(X,\mu)$ induces a Poisson suspension $\P(T)=T_\circ$, which is an automorphism of the space $(X_\circ,\mu_\circ)$:
$$\P(T)(x_1, x_2,\dots)=(Tx_1, Tx_2,\dots).$$
If the automorphism $T$, as a unitary operator acting in the space $L_2(X,\mu)$, has the property $T^n\to_w 0$, then the suspension $\P(T)$ has the mixing property:
$$ \mu_\circ(\P(T)^n C\cap C')\ \to \ \mu_\circ(C)\mu_\circ(C').$$

\section { Asymmetry of Automorphisms.}
The topic of asymmetry of automorphisms of Lebesgue spaces has a long history. Operators induced by an automorphism and its inverse are spectrally isomorphic. Automorphisms that are metrically (from the word "measure") not isomorphic to their inverse were discovered by Anzai \cite{An}, then pointed out in various classes of actions by Oseledets \cite{Os}, Ornstein and Shields \cite{Or}, and others.
Note that every asymmetric K-automorphism (examples are constructed in \cite{Or}) is not Bernoulli. We present some facts about asymmetry for non-mixing automorphisms. 

\newpage
\bf Theorem 1 (\cite{23}). \it The  generic automorphism $T$ of the  probability space has the following property: there exists a sequence
$m_j\to\infty$ such that for any measurable set $A$ the following convergences hold:
$$4\mu(A\cap T^{m_j}A\cap T^{3m_j}A)\to \ \mu(A)+ \mu(A)^2 + 2\mu(A)^3,$$
$$\mu(A\cap T^{-m_j}A\cap T^{-3m_j}A)\to \ \mu(A)^2.$$\rm

\vspace{3mm}
For infinite automorphisms, we have the following.

\vspace{3mm}
\bf Theorem 2 (\cite{R24}). \it In the group $Aut_\infty$, a dense $G_\delta$-set is formed by asymmetric automorphisms $T$ such that there exists a sequence
$m_j\to\infty$ such that for any set $A$ of finite measure
we have
$$3\mu(A\cap T^{m_j}A\cap T^{3m_j}A)\to \ \mu(A), \ \ \
\mu(A\cap T^{-m_j}A\cap T^{-3m_j}A)\to \ 0.$$\rm

\vspace{3mm}
To prove the genericity of this   property  we can show that all automorphisms satisfying the property form a $G_\delta$-set, find a suitable example of an automorphism, and use the density of its conjugacy class (see \cite{23}, \cite{R25}). Then  we use the fact that for a typical infinite automorphism $T$, there exists a sequence $n_j\to\infty$ such that
$T^{n_j}\to_w I/3$ (see Theorem 2.1 \cite{R25}). The presence of this weak limit implies the mutual singularity of  spectral measure $\sigma_T$  with all its convolution powers $\sigma_T^{\ast n}$ (see \cite{St}).

\bf Proposition (\cite{Ro}).  \it  If automorphisms $S$ and $T$ have spectral measures that are mutually singular with their convolution powers, respectively, then non-isomorphism of $S$ and $T$ implies non-isomorphism of the  suspensions $\P(S)$ and $\P(T)$.\rm 

Therefore, if the measure $\sigma_T$ is mutually singular with $\sigma_T^{\ast n}$ for all $n>1$, then non-isomorphism of $T$ and $T^{-1}$ implies non-isomorphism of $\P(T)$ and $\P(T^{-1})$. \rm

By virtue of what has been said, Poisson suspensions $\P(T)$ over a typical
automorphism $T$ (which, note, does not have the mixing property)
inherit the asymmetry property.
 
\newpage
\section{Mixing asymmetric Sidon automorphisms and Poisson suspensions}
Recall that the weak closure of a mixing action differs from the action itself by only one  point in operator space. For infinite measure space this  point is the zero operator, and for a probability space, it is the orthogonal projection onto the space of constant functions. We recall that the  mixing for infinite automorphism $T$ means  $T^n\to_w 0$ as $n\to\infty$.

\vspace{3mm}
\bf Theorem 3. \it There exists a mixing infinite automorphism $T$  and sequences
$m(j,i),\, 1\leq i\leq r_j$, $r_j\to\infty$ such that for any set $A$ of finite measure
$$\mu\left(\bigcup_{i=1}^{r_j-2}(A\cap T^{-m(j,i)}A\cap T^{m(j,i+1)}A)\right)\to \ 0,\eqno (0) $$
$$\mu\left(\bigcup_{i=1}^{r_j-2}(A\cap T^{m(j,i)}A\cap T^{-m(j,i+1)}A)\right)\to \ \mu(A).\eqno (1) $$
This automorphism $T$ is asymmetric.

\rm

\vspace{3mm}
Proof. Constructions of such automorphisms  can be found in \cite{24}. Recall the  definition of a rank-one automorphism  $T$. Its  parameters are
$h_1=1$, $r_j\geq 2$, $$ \bar s_j=(s_j(1), s_j(2),\dots,s_j(r_j)), \ s_j(i)\geq 0, \ j\in\N. $$
Its phase space $X$ is defined as  will be  the union of towers $$X_j=\bigcup_{i=0}^{h_j-1} T^iE_j,$$
where $T^iE_j$ are disjoint half-intervals called floors.

At stage $j$, the transformation $T$ is defined as  usual  translation of intervals, but at the top floor
$T^{h_j-1}E_j$, it is not yet defined.

Stage $j+1$.
Tower $X_j$ is cut into $r_j$ identical narrow subtowers $X_{j,i}$ (called columns), and $s_j(i)$ new floors are added above each column $X_{j,i}$. The transformation $ T $ is still defined as moving up one floor, but the last floor added above column number $ i $ is sent to the bottom floor of column number $ i + 1 $ . Thus, a new tower is created:
$$ X_{j + 1} = \bigcup_{i = 0} ^ { h_{j + 1} - 1} T^ i E_{j + 1}, \ \
h_{j + 1} = r_j h_j + \sum_{i = 1} ^ { r_j} s_j (i), $$
where $ E_{j + 1} $ is the bottom floor of column $ X_{j, 1} $ .
Note that by redefining the transformation, we completely preserve the previous constructions. After going through all the construction steps, we obtain an invertible transformation $ T: X \to X $ , which preserves the Lebesgue measure on the union $ X = \cup_j X_j $.

Let the construction $T$  the intersection
$X_j\cap T^mX_j$ for $h_{j}<m\leq h_{j+1}$ can be contained
in only one of the columns $X_{i,j}$ of the tower $X_j$. Such transformations
are called Sidon. The measure of  $X$ in this case is infinite.
As $r_j\to\infty$, then for  Sidon constructions and all sets $A,B$ of finite measure, we have
$$\mu(T^nA\cap B)\to 0, \ \ n\to\infty.$$
The Poisson suspension $\P(T)$  will 
be mixing with respect to the Poisson measure.

\bf Choice of parameters. \rm Let us  consider a construction with parameters $h_1=1$, $r_j=2^j$,
$$ \bar s_j=(dh_j, d^2h_j,\dots,d^{r_j}h_j), \ d>10, $$
setting $m(j,i)=h_j+s_j(i)$.
Let $A_j$ consist of the union of the floors of tower $X_j$.
Snce the intersection of $ A_j\cap T^{m(j,i)}A_j$ and the intersection of $ A_j\cap T^{-m(j,i+1)}A_j$ is $X_{j,i+1}\cap A_j$, for all $i=1,2,\dots,r_j-2$, we obtain
$$\mu\left(A_j\Delta \left(\bigcup_{i=1}^{r_j-2} (A_j\cap T^{m(j,i)}A_j\cap
T^{-m(j,i+1)}A_j)\right) \right) = 2 \mu(A_j)/ {r_j}. \eqno (2)$$

Any  set $A$ of finite measure is approximated
at stage $j$ by sets $A_j$ consisting of the union of some floors of tower $X_j$: $\mu(A\Delta A_j)\to 0$. Thus we get
$$\sum_{i=1}^{r_j} \mu\left((A\cap X_{j,i})\Delta (A_j\cap X_{j,i})\right)\ \to 0, \ j\to\infty.$$
For large $j$, for most $i$ ($1\leq i\leq r_j$), intersections of the form $A\cap T^{m(j,i)}A\cap T^{-m(j,i+1)}A$ differ little
from intersections $A_j\cap T^{m(j,i)}A_j\cap T^{-m(j,i+1)}A_j$.
For this reason, (2) implies (1).

Let us prove (0). Fix $j$. Let $A=C\uu D$, where $D=A\setminus X_j$.
For our Sidon construction, the equalities $$ X_j\cap T^{-m(j,i)}X_j =X_{j,i}, \ \ X_j\cap T^{m(j,i+1)}X_j =X_{j,i+2}, $$
but $X_{j,i+2}$ does not intersect $X_{j,i}$ ($i<r_j-1$),
therefore
$$\bigcup_{i=1}^{r_j-2} (X_j\cap T^{-m(j,i)}X_j\cap T^{m(j,i+1)}X_j ) = \emptyset. $$
Since $C\subset X_j$, we obtain
$$\bigcup_{i=1}^{r_j-2}(C\cap T^{-m(j,i)}C\cap T^{m(j,i+1)}C)=\emptyset. $$
Therefore, the set
$$Q_j:=\ \bigcup_{i=1}^{r_j-2}((C\uu D)\cap T^{-m(j,i)}(C\uu D)\cap T^{m(j,i+1)}(C\uu D)$$
included in the association
$$ D\cup \bigcup_{i=1}^{r_j-2}(C\cap T^{-m(j,i)}D ) \cup
\bigcup_{i=1}^{r_j-2}(C\cap T^{m(j,i+1)}D). $$
Let's show that
$$\mu\left( \bigcup_{i=1}^{r_j-2}(C\cap T^{-m(j,i)}D \right)\leq\mu(D). $$
To do this, we use the fact that all sets $ T^{m(j,i)}X_j$, $i=1,2, \dots, r_j-2$, are disjoint.
Let $D_i=D\cap T^{m(j,i)}X_j$, then, recalling that $C\subset X_j$, we obtain
$$\mu\left(\bigcup_{i=1}^{r_j-2}(C\cap T^{-m(j,i)}D)\right) \ = \
\sum_{i=1}^{r_j-2}\mu(T^{-m(j,i)}D_i)\ \leq\ \mu(D).$$
Similarly, we obtain
$$\mu\left( \bigcup_{i=1}^{r_j-2}(C\cap T^{m(j,i+1)}D\right)\leq\mu(D), $$
thus we obtained the inequality $\mu(Q_j)\leq 3 \mu(D)$. As $j$ grows, the measure $D=A\setminus X_j$ tends to 0, thereby establishing property (0). 

From equality
$$\mu\left(\bigcup_{i=1}^{r_j-2}(A\cap T^{-m(j,i)}A\cap T^{m(j,i+1)}A)\right) \ \ = $$ $$
\mu\left(\bigcup_{i=1}^{r_j-2}(SA\cap (STS^{-1})^{-m(j,i)}SA\cap (STS^{-1})^{m(j,i+1)}SA)\right)$$ 
it follows that The adjoint automorphism $STS^{-1}$ also has properties (0) and (1). It remains to note that, by (1), the automorphism $T^{-1}$ does not have property (0), so $T$ and $T^{-1}$ are not adjoint.
The theorem is proved.

\vspace{3mm}
\bf Theorem 4. \it Mixing Poisson suspensions for  some Sidon automorphisms are asymmetric. \rm

\vspace{3mm}
Proof.  The spectral measure of the Sidon construction described above is singular, and the convolution powers of this measure are absolutely continuous, as shown in \cite{24}. A sufficient condition for this spectral property to be realized is the inequality $\sum_j \frac 1 {r_j}<\infty$, which holds in our case since $r_j=2^j$. By virtue of Proposition 5.2 \cite{Ro} mentioned earlier, the asymmetry of the construction $T$ implies the asymmetry of its Poisson suspension $\P(T)$.

%\vspace{5mm}
%The author thanks the reviewer for their comments.
\large


\begin{thebibliography}{99}
\bibitem{KSF} Kornfeld I.P., Sinai Ya. G., Fomin S.V., Ergodic Theory, Nauka, Moscow, 1980

\bibitem{Ne}Neretin Yu.A., Categories of Symmetries and Infinite-Dimensional Groups, URSS, Moscow, 1998

\bibitem{PR} Parreau F., Roy E., Prime Poisson Suspensions, Ergodic Theory Dynam. Systems, 35:7 (2015), 2216-2230

\bibitem{R25} Ryzhikov V.V., Rigid Poisson Suspensions Without Roots, Mat. Zametki, 117:1 (2025), 146-150

\bibitem{An} Anzai, H., Ergodic skew product transformations on the torus, Osaka Math. J. 1951, Vol. 3, No. 1, pp. 83-99

\bibitem{Os} Oseledets, V. I., An example of two nonisomorphic systems with the same simple singular spectrum, Funct. Analysis and its Applications, 5:3 (1971), 75-79

\bibitem{Or} Ornstein, D. S., Shields, P. C., An uncountable family of K-automorphisms. Adv. Math. 10 (1973), 103-120

\bibitem{23} Ryzhikov, V. V.,  Generic extensions of ergodic systems, Sb. Math., 214:10 (2023), 1442-1457 

\bibitem{R24} Ryzhikov V.V., Generic properties of ergodic automorphisms, Tr. MMO, 85:2 (2024), 59-69

\bibitem{St} Stepin A.M., Spectral properties of typical dynamical systems, Izv. Akad. Sib. Mat., 50:4 (1986), 801-834

\bibitem{Ro} Roy E., Poisson suspensions and infinite ergodic theory, Ergodic Theory Dynam. Systems, 29:2 (2009), 667-683

\bibitem{24} Ryzhikov V.V., Polynomial rigidity and the spectrum of Sidon automorphisms. Mat. Sat., 215:7 (2024), 138-152

\end{thebibliography}
\end{document}